\documentclass{article}
\RequirePackage[amsthm,amsmath,natbib]{imsart}
\RequirePackage{imsart}
\RequirePackage[OT1]{fontenc}
\usepackage{amsmath,amssymb}

\usepackage[T1]{fontenc}
\usepackage{dsfont}
\begin{document}

\newcommand{\LL}{\mathcal{L}}
\newcommand{\und}{\vec{1}_d}
\newcommand{\suite}{_{n\geq 1}}
\newcommand{\ls}{\limsup}
\newcommand{\li}{\liminf}
\newcommand{\e}{\epsilon}
\newcommand{\eni}{\po \e_{n,i}\pf_{n\geq1, i \le p_n}}
\newcommand{\tab}{_{n\geq1,\;i\le p_n}}
\newcommand{\lb}{\newline}
\newcommand{\dd}{\delta}
\newcommand{\ld}{\llcorner}
\newcommand{\lr}{\lrcorner}
\newcommand{\alp}{\alpha}
\newcommand{\lab}{\lambda}
\newcommand{\Lab}{\Lambda}
\newcommand{\gam}{\gamma}
\newcommand{\Gam}{\Gamma}
\newcommand{\sig}{\sigma}
\newcommand{\Sig}{\Sigma}
\newcommand{\mmi}{\mid\mid}
\newcommand{\mmmi}{\mid\mid\mid}
\newcommand{\mmii}{\mid\mid_{[0,1]}}
\newcommand{\gMid}{\Bigg |}
\newcommand{\Mid}{\Big |}
\newcommand{\Mmi}{\Big|\Big|}
\newcommand{\Mmii}{\Big|\Big|_{[0,1]}}
\newcommand{\T}{\Theta}
\newcommand{\tht}{\theta}
\newcommand{\sli}{\sum\limits}
\newcommand{\sliin}{\sum\limits_{i=1}^n}
\newcommand{\sliid}{\sum\limits_{i=1}^d}
\newcommand{\sliik}{\sum\limits_{i=1}^k}
\newcommand{\sliiN}{\sum\limits_{i=1}^N}
\newcommand{\slijk}{\sum\limits_{j=1}^k}
\newcommand{\proliin}{\prod\limits_{i=1}^n}
\newcommand{\proliid}{\prod\limits_{i=1}^d}
\newcommand{\proliik}{\prod\limits_{i=1}^k}
\newcommand{\proliiN}{\prod\limits_{i=1}^N}
\newcommand{\prolijk}{\prod\limits_{j=1}^k}
\newcommand{\ili}{\int\limits}
\newcommand{\proli}{\prod\limits}
\newcommand{\limn}{\lim_{n\rightarrow\infty}\;}
\newcommand{\limk}{\lim_{k\rightarrow\infty}\;}
\newcommand{\inv}{\frac{1}}
\newcommand{\lsn}{\limsup_{n\rightarrow\infty}}
\newcommand{\lin}{\liminf_{n\rightarrow\infty}}
\newcommand{\lsk}{\limsup_{k\rightarrow\infty}}
\newcommand{\lik}{\liminf_{k\rightarrow\infty}}
\newcommand{\bculi}{\bigcup\limits}
\newcommand{\bcali}{\bigcap\limits}
\newcommand{\wt}{\widetilde}
\newcommand{\indep}{{\bot}\kern-0.9em{\bot}}
\newcommand{\beq}{\begin{equation} }

\newcommand{\eeq}{\end{equation} }
\newcommand{\mcal}{\mathcal}
\newcommand{\sq}{\sqrt}
\newcommand{\rar}{\rightarrow}
\newcommand{\lar}{\leftarrow}
\newcommand{\cvps}{\rightarrow_{p.s.}\;}
\newcommand{\cvpr}{\rightarrow_{P}\;}
\newcommand{\cvloi}{\rightarrow_{\mcal{L}}\;}
\newcommand{\ii}{\vec{i}}
\newcommand{\jj}{\vec{j}}
\newcommand{\kk}{\vec{k}}
\newcommand{\zzi}{z_{\vec{i}}}
\newcommand{\zzij}{z_{\vec{i},\vec{j}}}
\newcommand{\zzI}{\vec{z}_I}
\newcommand{\zzIJ}{\vec{z}_{I,J}}
\newcommand{\zz}{z}
\newcommand{\nk}{n_k}
\newcommand{\nkm}{n_{k-1}}
\newcommand{\nkk}{n_{k+1}}
\newcommand{\nkd}{{n_k}}
\newcommand{\nkkd}{{n_{k+1}}}
\newcommand{\Ln}{L_{n,1}}
\newcommand{\Lnn}{L_{n,2}}
\newcommand{\Lnl}{L_{n,\ell}}
\newcommand{\Pn}{\Pi_{n,1}}
\newcommand{\Pnn}{\Pi_{n,2}}
\newcommand{\Pnl}{\Pi_{n,\ell}}
\newcommand{\mn}{\vec{m}_n}
\newcommand{\Mn}{\vec{M}_n}
\newcommand{\hhn}{h_{n,1}}
\newcommand{\hhnn}{h_{n,2}}
\newcommand{\hhk}{h_{n_k,1}}
\newcommand{\hhkk}{h_{n_k,2}}
\newcommand{\mk}{\vec{m}_k}
\newcommand{\Mk}{\vec{M}_k}
\newcommand{\hhnl}{h_{n,\ell}}
\newcommand{\hnk}{h_{n_k}}
\newcommand{\bn}{b_{n,1}}
\newcommand{\bnn}{b_{n,2}}
\newcommand{\bnl}{b_{n,\ell}}
\newcommand{\bnk}{b_{n_k}}
\newcommand{\bnkk}{b_{n_{k+1}}}
\newcommand{\hkl}{{h_{n_k}^{(l)}}}
\newcommand{\hnl}{{h_n^{(l)}}}
\newcommand{\Tn}{T_{n,1}}
\newcommand{\Tnn}{T_{n,2}}
\newcommand{\Tnl}{T_{n,\ll}}
\newcommand{\hamn}{\hat{m}_n}
\newcommand{\harn}{\hat{r}_n}
\newcommand{\hafn}{\hat{f}_n}
\newcommand{\srl}{\stackrel}
\newcommand{\wap}{(\Omega,\mathcal{A},\rm I\kern-2pt P)}
\newcommand{\aoo}{\Big\{}
\newcommand{\aff}{\Big\}}
\newcommand{\coo}{\Big [}
\newcommand{\cff}{\Big]}
\newcommand{\poo}{\Big (}
\newcommand{\pff}{\Big)}
\newcommand{\po}{\big (}
\newcommand{\pf}{\big)}
\newcommand{\ao}{\big \{}
\newcommand{\af}{\big \}}
\newcommand{\co}{\big [}
\newcommand{\cf}{\big ]}
\newcommand{\pooo}{\bigg (}
\newcommand{\pfff}{\bigg)}
\newcommand{\aooo}{\bigg \{}
\newcommand{\afff}{\bigg \}}
\newcommand{\cooo}{\bigg [}
\newcommand{\cfff}{\bigg ]}
\newcommand{\poooo}{\Bigg (}
\newcommand{\pffff}{\Bigg)}
\newcommand{\aoooo}{\Bigg \{}
\newcommand{\affff}{\Bigg \}}
\newcommand{\coooo}{\Bigg [}
\newcommand{\cffff}{\Bigg ]}
\newcommand{\FF}{\mathcal{F}}
\newcommand{\TT}{\mathcal{T}}
\newcommand{\GG}{\mathcal{G}}
\newcommand{\BBGG}{\mathcal{B}(\mathcal{G})}
\newcommand{\CC}{\mathcal{C}}
\newcommand{\KK}{\mathcal{K}}
\newcommand{\SSS}{\mathcal{S}}
\newcommand{\BB}{\mathcal{B}}
\newcommand{\PP}{\mathcal{P}}
\newcommand{\HH}{\mathcal{H}}
\newcommand{\NN}{\mathcal{N}}
\newcommand{\MM}{\mathcal{M}}
\newcommand{\DD}{\mathcal{D}}
\newcommand{\cc}{\widetilde{c}}
\newcommand{\EEE}{\mathbb{E}}
\newcommand{\NNN}{\mathbb{N}}
\newcommand{\PPP}{ \mathbb{P}}
\newcommand{\CCC}{\mathbb{C}}
\newcommand{\KKK}{\mathbb{K}}
\newcommand{\RRR}{\mathbb{R}}
\newcommand{\wtc}{\widetilde{c}}
\newcommand{\wtf}{\widetilde{f}}
\newcommand{\wth}{\widetilde{h}}
\newcommand{\wtn}{\widetilde{n}}
\newcommand{\wtv}{\widetilde{v}}
\newcommand{\wtA}{\widetilde{A}}
\newcommand{\wtC}{\widetilde{C}}
\newcommand{\wtE}{\widetilde{E}}
\newcommand{\ovg}{\overline{g}}
\newcommand{\ovh}{\overline{h}}
\newcommand{\ovE}{\overline{E}}
\newcommand{\wtG}{\widetilde{G}}
\newcommand{\ovH}{\overline{H}}
\newcommand{\ovI}{\overline{I}}
\newcommand{\ovJ}{\overline{J}}
\newcommand{\ovK}{\overline{K}}
\newcommand{\wtN}{\widetilde{N}}
\newcommand{\wtP}{\widetilde{P}}
\newcommand{\ovR}{\overline{R}}
\newcommand{\wtPPP}{\widetilde{\mathbb{P}}}
\newcommand{\ovPPP}{\overline{\mathbb{P}}}
\newcommand{\wtF}{\widetilde{F}}
\newcommand{\wtI}{\widetilde{I}}
\newcommand{\wtK}{\widetilde{K}}
\newcommand{\wtFF}{\widetilde{\mathcal{F}}}
\newcommand{\ovFF}{\overline{\mathcal{F}}}
\newcommand{\cov}{\mathrm{Cov}}
\newcommand{\kif}{k\rightarrow\infty}
\newcommand{\nif}{n\rightarrow\infty}
\newcommand{\FFGG}{\FF\times\GG}
\newcommand{\vk}{\vskip10pt}
\newcommand{\Tproj}{\mathcal{T}_{0}}
\newcommand{\TTd}{\mathcal{T}_d}
\newcommand{\farc}{\frac}
\newcommand{\Nf}{\nabla_f}
\newcommand{\Nfn}{\nabla_f(\log_2n )}
\newcommand{\Cf}{\chi_f}
\newcommand{\nono}{\nonumber}
\newcommand{\Dn}{\Delta_n}
\newcommand{\DPn}{\Delta\Pi_n}
\newcommand{\norm}{\mid\mid \cdot \mid\mid}
\newcommand{\Xt}{(X(t))_{t\in T}}
\newcommand{\zklj}{z_{k,l,j}}
\newcommand{\znlj}{z_{n,l,j}}
\newcommand{\Idd}{[0,1[^d}
\newcommand{\iideux}{\vec{i}\in {\{1,\ldots,2^p\}}^d}
\newcommand{\hk}{h_{n_k}}
\newcommand{\zik}{z_{i,n_{k}}}
\newtheorem{theo}{Theorem}
\newtheorem{ptheo}{Proof of theorem}
\newtheorem{lem}{Lemma}[section]
\newtheorem{plem}{Proof of lemma}[section]
\newtheorem{prop}{Property}[section]
\newtheorem{preuveprop}{Proof of property}[section]
\newtheorem{defi}{Definition}[section]
\newtheorem{propo}{Proposition}[section]
\newtheorem{popo}{Proof of proposition}[section]
\newtheorem{coro}{Corollary}[theo]
\newtheorem{pcoro}{Proof of corollary}[theo]
\newtheorem{rem}{Remark}[subsection]
\newtheorem{Fact}{Fact}[section] \numberwithin{equation}{section}
\newtheorem{ineg}{Inequality}[section]
\begin{frontmatter}
\title{Clustering rates and Chung type functional laws of the iterated logarithm for empirical and quantile processes}
\runtitle{Clustering rates and Chung limit laws for empirical processes}
\affiliation{}

\author{\fnms{Davit} \snm{Varron}\ead[label=e1]{davit.varron@univ-fcomte.fr}}
\address{Laboratoire de Mathématiques de Besançon, UMR CNRS 6623,\break Université de Franche-Comté\break \printead{e1}}

\runauthor{D. Varron}

\maketitle

\begin{abstract}
Following the works of Berthet \cite{Berthet1,Berthet2}, we first
obtain exact clustering rates in the functional law of the
iterated logarithm for the uniform empirical and quantile processes and for their
increments. In a second time, we obtain
functional Chung-type limit laws for the local empirical process
for a class of target functions on the border of the Strassen set.
\end{abstract}


\end{frontmatter}
\section{Introduction} Define the
uniform empirical process by $\alp_n(t):=n^{1/2}(F_n(t)-t)$, where
$F_n(t):=n^{-1}\sharp\ao i\in \{1,\ldots,n\},\;U_i\le t\af,\;t\in [0,1]$, and
\((U_i)\suite\) are independent, identically distributed (i.i.d) random variables uniformly distributed
on $[0,1]$. Define the quantile process by
\[\beta_n(t)=n^{1/2}\Big(F_n^{-1}(t)-t\Big),\;t\in [0,1],\]
where \(F_n^{-1}(t):=\inf\{u:\;F_n(u)\geq t\}\). In a
metric space $(\mcal{E},d)$ we write $u_n\leadsto \mcal{H}$ whenever
$u_n$ is relatively compact with limit set $\mcal{H}$ (see, e.g.,
\cite{Mason2}). The two above mentioned processes
have been extensively investigated in the literature (see, e.g.,
\cite{MSW} and \cite{Vander} and the references therein). In a
pioneering work, Finkelstein \cite{Finkel} has established the
functional law of the iterated logarithm (FLIL) for $\alp_n$.
Namely, the author showed that, writing \(\log_2u=\log(\log (u\vee e))\) and $ b_n=\sqrt{2\log_2n }$, we have~:
\beq\frac{\alp_n}{b_n}\leadsto_{a.s.}\;\SSS_2\label{cmua},  \eeq
in the metric space $(B[0,1],\mmi
\cdot \mmi)$, where $B[0,1]$ stands for the set of bounded functions
on $[0,1]$ and $\mmi \cdot \mmi$ is the sup-norm over $[0,1]$. The
set $\SSS_2$ in (\ref{cmua}) is given by \beq \SSS_2:=\Big \{f(t)\in
\SSS_1,\;f(1)=0\Big \}\label{sss2},\eeq where \beq \SSS_1:=\aooo
f\in B[0,1],\; \exists f' \text{ Borel, }f(\cdot):=\ili_0^\cdot f'(t)dt,\;\ili_0^1f'^2(t)dt\le 1\afff \label{sss1}.\eeq
Note that $\SSS_2$ (resp. $\SSS_1$) is the unit ball of the
reproducing kernel Hilbert space of the Brownian bridge (resp. of
the Wiener process) on [0,1]. In the spirit of \cite{Finkel}, Mason
\cite{Mason2} has obtained the following FLIL for the local
empirical process~: \beq\frac{\alp_n(a_n\cdot
)}{\sqrt{a_n}b_n}\leadsto_{a.s.} \SSS_1\label{flilmason}.\eeq Here,
\(a_n\) is a sequence of constants satisfying $a_n \downarrow
0,\;na_n\uparrow \infty$ and $na_n/\log_2n\rar\infty$. Deheuvels and
Mason \cite{DeheuvelsM2} have established a related uniform
functional limit law for the following collections of random trajectories.
\[\T_n:=\Big \{\frac{\alp_n(t+a_n\cdot )-\alp_n(t)}{\sqrt{2a_n\log(1/a_n)}},\;t\in [0,1-a_n]\Big \}.\]
They showed that, with probability one :
\[\limn \sup_{g_n\in \T_n}\inf_{f\in \SSS_1} \mmi g_n-f\mmi=0,\]
\beq \limn \sup_{f\in \SSS_1}\inf_{g_n\in \T_n} \mmi
g_n-f\mmi=0\label{dm2},\eeq where $a_n$ is a sequence of constants
fulfilling $a_n \downarrow 0,\;na_n\uparrow \infty,\; na_n/\log
n\rar\infty,\;\log(1/a_n)/\log_2n\rar \infty$. Berthet
\cite{Berthet1} refined (\ref{dm2}) under slightly stronger conditions
imposed upon \(a_n\). Making use of sharp upper bounds for
Gaussian measures due to Talagrand \cite{Tal1}, he proved that
for any \(\e>\e_0\) (where $ \e_0 $ is a universal
constant), we have almost surely for all $n$ large enough~: 
\beq\T_n\subset \SSS_1+\e
\log(1/a_n)^{-2/3}\mcal{B}_0\label{berthet1}.\eeq Here
\(\mcal{B}_0:=\{f\in B[0,1]: \;\mmi f\mmi\le 1\}\). The first aim of
the present article is to show that the techniques employed in the just-mentioned result can be adapted to some other random objects than that
used for that given in (\ref{berthet1}) (see Theorems \ref{theo:
recfinkel} and \ref{theo: recproclocal} in the sequel). Results of
this kind are usually called \textit{clustering rates}. Another related problem
is finding rates of convergence of such random sequences to a
specified function belonging to the cluster set. Such results are
known under the name of functional \textit{Chung-type limit laws}.
We now focus on the local empirical process $\alp_n(a_n\cdot)$,
where $a_n\downarrow 0$ as $\nif$. The works of Cs\'aki
\cite{Csaki}, de Acosta \cite{Acosta1}, Grill \cite{Grill1}, Gorn
and Lifshits \cite{GornLif}, and Berthet and Lifshits
\cite{BerthetLif} on small ball probabilities for Wiener processes
provide some crucial tools to establish such limit laws for
$(\alp_n(a_n\cdot))\suite$, as these are expected to
asymptotically mimic their gaussian analogues (see Mason \cite{Mason2}).
Along this line, Deheuvels \cite{Deheuvels2} established Chung-type
limit laws for $(\alp_n(a_n\cdot))\suite$, by showing that, if
$a_n$ is a sequence of constants satisfying $na_n \uparrow
\infty,\;a_n \downarrow 0$ and $na_n/(\log_2 n)^3\rar\infty$, we
have, almost surely, for each $f\in \SSS_1$ satisfying $\mmi f
\mmi_H^2:=\ili_0^1 f'^2(t)dt<1$~:
\[\lin\; (\log_2n) \Mmi \frac{\alp_n(a_n
\cdot)}{\sqrt{a_n}b_n}-f\Mmi=\frac{\pi}{4\sqrt{1-\mmi f\mmi_{H}}}.\]
The proof of this theorem relies on strong approximation methods in
combination with the results of de Acosta \cite{Acosta1}. The latter
provides useful exponential bounds for
$$\PPP\poo \Mmi \frac{W}{T}-f\Mmi\le \e\pff,$$ with a small $\e>0$ and a large
$T$. Here, $W$ is a Wiener process on $[0,1]$ and $f$ satisfies
$\mmi f \mmi_{H}^2<1$. The study of related probabilities when $\mmi
f \mmi_{H}=1$ has required different arguments. In \cite{Grill1},
rough estimates are given. In \cite{GornLif} and
\cite{BerthetLif}, some exact rates are given, but only for functions with
first derivatives having a variation either bounded or locally
infinite. The sets of all functions of this type are called
 $\SSS_1^{bv}$ and $\SSS_1^{liv}$ respectively. In the present paper,
we shall make use of the latter results to extend the work of
Deheuvels \cite{Deheuvels2} to the case where $f\in \SSS_1^{bv}\cup\SSS_1^{liv}$. The remainder of our paper is
organized as follows. Our main results are stated in $\S 2$,
Theorems \ref{theo: recfinkel}, \ref{theo: recproclocal} and
\ref{theo: Chunglocal}. In $\S 3$, the proofs of these theorems
are provided.
\section{Main Results}
Our first result gives clustering rates in Finkelstein's FLIL
\cite{Finkel}.
\begin{theo}\label{theo: recfinkel}
There exists a universal constant \(\e_0\) >0 such that, for
any
choice of \(\e>\e_0\) we have almost surely, for all large $n$\begin{align}\frac{\alp_n}{(2\log_2n )^{1/2}} &\in\;\SSS_2+\e (\log_2n) ^{-2/3} \mcal{B}_0\label{recfinkel},\\
\frac{\beta_n}{(2\log_2n )^{1/2}} &\in\; \SSS_2+\e (\log_2n) ^{-2/3}
\mcal{B}_0\label{recquantile}.
\end{align}

\end{theo}
\begin{rem}
The uniform Bahadur-Kiefer representation (see \cite{BaKief}) asserts that, almost surely~: 
\[\lsn n^{1/4}(\log n)^{-1/2}(\log_2 n)^{-1/4}\mmi \alp_n+\beta_n\mmi=2^{-1/4},\]
from where (\ref{recquantile}) is readily implied by
(\ref{recfinkel}).\end{rem} Our second theorem concerns the
FLIL for local increments of the empirical process.
\begin{theo}\label{theo: recproclocal}
Let \(a_n\) be positive real numbers satisfying, as $n\rar
\infty$, \beq na_n\uparrow\infty,\;\;\frac{na_n}{(\log_2n)
^{7/3}}\rar\infty,\;\;a_n\downarrow 0\label{sirn}.\eeq Then there
exists a universal constant \(\e_1>0\) such that, for any choice of
\(\e> \e_1\) we have almost surely, for all large $n$, \beq
\frac{\alp_n(a_n\cdot )}{\sqrt{2 a_n \log_2n }}\in \SSS_1+\e
(\log_2n) ^{-2/3}\mcal{B}_0\label{recproclocal}.\eeq If moreover
$na_n/(\log_2 n)^{11/3}\rar \infty$ then we have, almost surely,
ultimately as $\nif$, \beq\frac{\beta_n(a_n\cdot )}{\sqrt{2 a_n
\log_2n }}\in\SSS_1+\e (\log_2n)
^{-2/3}\mcal{B}_0\label{recquantilelocal}.\eeq
\end{theo}
\begin{rem}
We shall use the fact (see e.g. \cite{EinmahlM2}, Theorem 5) that, under
(\ref{sirn}), we have almost surely
\[\lsn \;(n/a_n)^{1/4}(\log_2 n)^{-1/4}(2\log_2 n+\log(na_n))^{-1/2}\mmi \alp_n(a_n\cdot)+\beta_n(a_n\cdot)\mmi\le 2^{-1/4},\]
from where  (\ref{recquantilelocal}) is implied by
(\ref{recproclocal}) after straightforward
computations.\end{rem} In order to state our last result, we
need to give some definitions. Recall that $f\in\SSS_1^{bv}$
whenever $f'$ has a derivative with bounded variation and
$\ili_0^1 {f'}^2(t)dt =1$. Results on small ball probabilities
for a Wiener process when $f\in \SSS_1^{bv}$ have been established
by Gorn and Lifshits \cite{GornLif}. For such a function $f$, we
shall write $\Nf(L):=L^{2/3}, \;L>0$ and we denote by $\Cf$ the constant
which is the unique solution of equation (3.1) in \cite{GornLif}
(we refer to the just mentioned paper for more details). The case
where $f\in \SSS_1^{liv}$ (i.e. where $\ili_0^1 {f'}^2(t)dt =1$ and the derivative of $f'$
admits a version with locally infinite variation) has been treated
by Berthet and Lifshits \cite{BerthetLif}. For such a function
$f$, we set $\Cf:=1$ and we denote by $\Nf(L)$ the unique solution
of equation (2.1) in \cite{Berthet2}. Our third result is stated as follows.
\begin{theo}\label{theo: Chunglocal} Let $f\in \SSS_1^{bv}\cup\SSS_1^{liv}$ be arbitrary and let
$a_n$ be a sequence of real numbers satisfying, as $n\rar\infty$,
\beq na_n \uparrow \infty,\;a_n\downarrow 0,\; a_n \log_2n \rar
0,\label{ww}\eeq \beq\limn \frac{na_n}{\log_2n \Nf^2(\log_2n
)}=\infty.\label{www}\eeq Then we have, almost surely~:  \[\lin
\Nf(\log_2n ) \Mmi\frac{\alp_n(a_n \cdot)}{\sqrt{2a_n\log_2n
}}-f\Mmi=\Cf.\]
\end{theo}
\begin{rem} The conditions (\ref{ww}) and (\ref{www}) imposed upon
 $a_n$ turn out to be the best possible with respect to the methods used in the proof of Theorem \ref{theo: Chunglocal}. The latter combines poissonization techniques with strong
 approximation arguments. Deheuvels and
Lifshits \cite{DeheuvelsLif3} and Shmileva \cite{Shmileva1} have
provided new tools to estimate probabilities of shifted small balls
for a Poisson process without making use of strong approximation
techniques. These results show up to be powerful enough to
investigate Chung-Mogulskii limit laws for $\alp_n(a_n.)$ without
making use of strong approximation techniques, and thus relaxing
condition (\ref{ww}). However, the just-mentioned results do not cover
the case where $f\in \SSS_1^{liv}$.
\end{rem}
\section{Proofs}
\subsection{Proof of Theorem \ref{theo: recfinkel}}\label{uol}
Select an \(\e>0\) and consider the sequence \(\e_n:= \e (\log_2n)
^{-2/3}\). The main tool to achieve our goal is the exponential
inequality stated in the following fact, which follows directly from
Talagrand \cite{Tal3}. Recall that $\mcal{B}_0$ is the unit ball for
$\norm$.
\begin{Fact}\label{propo: recouvrement Pbrownien}
Let $B$ be a Brownian bridge. There exists three constants $K$,
$L_0$ and \(u_0\)
>0 such that, for any \(0<u<u_0\) and $c>0$, we have~:  \beq P\po
B\notin c \SSS_2+u \mcal{B}_0\pf\le
K\exp\Big(\frac{L_0}{u^2}-\frac{c
u}{2}-\frac{c^2}{2}\Big).\label{rpb}\eeq Let $W$ be a Wiener process
on $[0,1]$. There exist two constants $u_1$ and $L_1$ such that, for
any \(0<u<u_1\) and $c>0$, we have \beq P\po W\notin c \SSS_1+u
\mcal{B}_0\pf\le \exp\Big(\frac{L_1}{u^2}-\frac{c
u}{2}-\frac{c^2}{2}\Big).\label{rw}\eeq
\end{Fact}
In the proof of Theorem \ref{theo: recfinkel}, we will make use of
blocking techniques (see, e.g., \cite{DeheuvelsM2} and
\cite{Berthet1}). For any real umber $a$, set $[a]$ as the unique
integer $q$ fulfilling $q\le a<q+1$, and set
$$\nkd:=\Big{[}\exp\poo k\exp\poo-(\log k)^{1/6}\pff\pff\Big{]},\;k\geq 1.$$
Set $N_k:= \{\nkd,\ldots ,\nkkd-1 \}$ for $k\geq 5$. Given an
integer $n\geq 1$, we set \(k(n)\) as the unique integer $k$ such
that \(n\in N_k\).  We shall first study the following sequence of
functions
\[g_n:=(\nkkd)^{-1/2} b_{\nkkd}^{-1}H_n,\;k=k(n),\]
with $H_n(t):=n (F_n(t)-t)$ and $b_n:=(2\log_2 n)^{1/2}$.
 Let \(p_1\) and \(q_1\) be two
conjugates numbers (such that \(1/p_1+1/q_1=1\) ) with \(1<p_1<\infty\).
Set, for $k\geq 1$,
\[m_{p_1,k}:= \min_{n\in N_k} \PPP\pooo\frac{1}{(\nkkd)^{1/2} b_{\nkkd}}\mmi H_{\nkkd}-H_n \mmi\le \frac{1}{p_1}\e_{\nkkd}\pfff.\]
A standard blocking argument based upon Ottaviani's inequality (see,
e.g., \cite{DeheuvelsM2}, Lemma 3.4) yields
\begin{align}
  \nonumber&\PPP\pooo\bculi_{n\in N_k}\Big \{\frac{1}{(\nkkd)^{1/2} b_{\nkkd}} H_n\notin \SSS_2+\e_{\nkkd}\mcal{B}_0\Big \}\pfff \\ \nonumber
  \nonumber\le \;& \frac{1}{m_{p_1,k}}\PPP\pooo\frac{1}{(\nkkd)^{1/2} b_{\nkkd}}
H_{\nkkd}\notin \SSS_2+\frac{1}{q_1}\e_{\nkkd}\mcal{B}_0\pfff.
\end{align}
Let $k$ be integer and select \(n\in N_k\). By the
Dvoretsky-Kiefer-Wolfowitz inequality (see, e.g., \cite{Tal2}) we
have~: 
\begin{align}
\nonumber&\PPP \pooo\frac{1}{(\nkkd)^{1/2} b_{\nkkd}} \mmi H_{\nkkd}-
H_n \mmi \geq \frac{1}{p_1}\e_{\nkkd}\pfff\\ \nonumber \nonumber\le\;
& \PPP\pooo\mmi \alp_{\nkkd-n}\mmi\geq
\frac{1}{p_1}\e_{\nkkd}\po\frac{1}{1-\frac{\nkd}{\nkkd}}\pf^{1/2}b_{\nkkd}\pfff\\
\nonumber \nonumber\le\; & 3\exp\poo-\frac{4\e^2}{p_1^2}\frac{\log_2(n_k)^{-1/3}}{\po 1-\frac{\nkd}{\nkkd}\pf}\pff \text{ for large enough } k,
\end{align}whence $m_{p_1,k}\geq 1/2$ for all large $k$ by routine analysis. Now let \(p_2,q_2>1\) be two conjugate numbers. For $k\geq 1$ we have, $B_{\nkkd}$ denoting a
Brownian bridge,
\begin{align}
\nonumber&\PPP\pooo\frac{\alp_{\nkkd}}{b_{\nkkd}}\notin
\SSS_2+\frac{1}{q_1}\e_{\nkkd} \mcal{B}_0\pfff \nonumber
\nonumber\le\; \PPP\pooo\mmi \alp_{\nkkd} -B_{\nkkd}\mmi \geq
\frac{1}{p_2q_1}\e_{\nkkd}b_{\nkkd}\pfff \\
\nonumber &\;+\PPP\pooo B_{\nkkd}\notin
b_{\nkkd}\SSS_2+\frac{1}{q_2q_1}\e_{\nkkd}b_{\nkkd} \mcal{B}_0\pfff\\
\nonumber \nonumber:=\;&\PPP_k^{KMT}+\PPP_k^{Tal}.
\end{align}
Making use of the Koml\'os-Major-Tusnàdy approximation (see, e.g., \cite{KMT})), we can
choose a sequence $(B_{\nkd})_{k\geq 1}$ satisfying, for some
universal constants $C_2,C_3$ and for all $k$ large enough,
\[\PPP_k^{KMT}\le C_2\exp\Big(-C_3(\nkkd)^{1/2}\frac{1}{2p_2q_1}\e_{\nkkd}b_{\nkkd}\Big).\]
On the other hand, by applying assertion  (\ref{rpb}) of Fact \ref{propo: recouvrement Pbrownien} we have, for all
large $k$,
\[\PPP_k^{Tal}\le K \exp\Big[-\Big(\frac{\e}{q_1q_2}-\frac{L_0(q_1q_2)^2}{2\e^2}\Big)(\log_2 \nkkd)^{1/3}-\log_2 \nkkd\Big].\]
Routine analysis shows that both $\PPP_k^{KMT}$ and $\PPP_k^{Tal}$
are sumable in $k$ for any choice of
$\e>(L_0/2)^{1/3}=:\e_0$, provided that $q_1,q_2$ are chosen close enough to 1. Now an application of (\ref{cmua}) in combination with
elementary properties of the sequence $(\nkd)_{k\geq 1}$ shows that, almost surely,
as $\nif$,
\[\mmi g_n-b_n^{-1}\alp_n \mmi=o((\log_2n) ^{-2/3}).\]
\subsection{Proof of Theorem \ref{theo: recproclocal}}
Recall that $b_n:= (2\log_2 n)^{1/2},\;n\geq1$. Let \(p_1,q_1>1\) be
two conjugate numbers. Set, for $k\geq 1$~: 
\[\mathfrak{m}_{p_1,k}:= \min_{n\in N_k}\PPP\poo\frac{1}{(\nkkd a_{\nkkd})^{1/2}b_{\nkkd}}\mmi H_n(a_{\nkkd}\cdot )-H_{\nkkd}(a_{\nkkd}\cdot) \mmi\le \frac{1}{p_1}\e_{\nkkd}\pff.\]
The same blocking argument as in \S \ref{uol} yields
\[\PPP\pooo\bculi_{n\in N_k}\Big \{\frac{H_n(a_{\nkkd}\cdot )}{(\nkkd a_{\nkkd})^{1/2}b_{\nkkd}}\notin \SSS_1+\e_{\nkkd} \BB_0\Big\}\pfff\]\[\le \frac{1}{\mathfrak{m}_{p_1,k}} \PPP\poo\frac{H_{\nkkd}(a_{\nkkd}\cdot )}{(\nkkd a_{\nkkd})^{1/2}b_{\nkkd}}\notin \SSS_1+\frac{1}{q_1}\e_{\nkkd}\BB_0\pff.\]
Now, for any integer $k\geq 5$ and $n\in N_k$, we have
\begin{align}
\nono &\PPP\poo\frac{1}{\sqrt{\nkkd}}\mmi H_{\nkkd}(a_{\nkkd}\cdot
)-H_n(a_{\nkkd}\cdot ) \mmi\geq
\frac{1}{p_1}\e_{\nkkd}b_{\nkkd}\pff\\
\nono\le &\PPP\poo\sup_{t\le a_{\nkkd}}\frac{\mid
\alp_{\nkkd-n}(t)\mid}{1-t}\geq
\frac{1}{p_1}\e_{\nkkd}b_{\nkkd}\po\frac{\nkkd
a_{\nkkd}}{\nkkd-\nkd}\pf^{1/2}\pff.\end{align}
 It is well known
(see, e.g., \cite{MSW}, Proposition 1, p. 133) that for each $n$,
the process \((1-t)^{-1}\alp_n(t)\) is a martingale in $t$. The
Doob-Kolmogorov inequality yields~:
\[1-\mathfrak{m}_{p_1,k}\le \frac{p_1^2(1-a_{\nkkd})(1-\frac{\nkd}{\nkkd})(\log \nkkd)^{1/3}}{2\e^2}.\]
Hence for all $k$ large enough we have $\mathfrak{m}_{p_1,k}\geq 1/2$. Now set
for each integer $n\geq
1$, \[\widetilde{\Pi}_n(t):=n^{-1/2}\poo\sli_{i=1}^{\eta_n}\mathds{1}_{\{U_i\le
t\}}-t\pff ,\;t\in[0,1],\] where \(\eta_n\) is a Poisson variable
with expectation $n$, which is independent of
 \((U_i)_{i\geq 1}\). Let $\wt{\Pi}$ denote a standard centered Poisson process on $\RRR^+$ and let $W$ be a Wiener process that we assume to be constructed on the same underlying probability space as $\wt{\Pi}$.
 Notice that $\wt{\Pi}_n(\cdot)$ and $n^{-1/2}\wt{\Pi}(n\cdot)$ are equal in distribution as processes on $[0,1]$.
 Now let $p_2,q_2>1$ be two conjugate numbers. By making use of
Poissonization techniques (see, e.g., \cite{DeheuvelsM2}, Lemma 2.1 or \cite{VarronSPA}, Proposition 2.1 for a more general form) we see
that, for all sufficiently large $k$~: 
\begin{align}
\nonumber&\PPP\pooo\frac{\alp_{\nkkd}(a_{\nkkd}\cdot
)}{a_{\nkkd}^{1/2}b_{\nkkd}}\notin \SSS_1+\frac{1}{q_1}\e_{\nkkd}
\mcal{B}_0\pfff\\ \nonumber \nonumber\le
\;&2\PPP\pooo\frac{\wt{\Pi}_{\nkkd}(a_{\nkkd}\cdot
)}{a_{\nkkd}^{1/2}b_{\nkkd}}\notin \SSS_1+\frac{1}{q_1}\e_{\nkkd}
\mcal{B}_0\pfff\\ \nonumber=\;&2\PPP\pooo\frac{\wt{\Pi}(\nkkd
a_{\nkkd}\cdot)}{(2\nkkd a_{\nkkd}\log_2
\nkkd)^{1/2}}\notin \SSS_1+\frac{1}{q_1}\e_{\nkkd} \mcal{B}_0\pfff\\
\nonumber \nonumber\le\; & 2\PPP\Big(\mmi W(\nkkd
a_{\nkkd}\cdot)-\wt{\Pi}(\nkkd a_{\nkkd} \cdot)\mmi \geq
\frac{1}{q_1p_2}(\nkkd a_{\nkkd})^{1/2}b_{\nkkd}\e_{\nkkd}\Big)\\
\nonumber \nonumber&+2\PPP\Big(\frac{W(\nkkd a_{\nkkd}\cdot)}{(\nkkd
a_{\nkkd})^{1/2}b_{\nkkd}}\notin \SSS_1+ \frac{1}{q_1
q_2}\e_{\nkkd}\Big)\\ \nonumber \nonumber:=\;&
\PPP_k^{KMT}+\PPP_k^{Tal}. \end{align} Now, making use of the strong
approximation theorem of Koml\'os-Major-Tusnàdy \cite{KMT2}, we
can assume that the process $W$ involved in the former expression
satisfies, for some universal constants $C_1,C_2,C_3>0$, and for all
$T>0,\;z>0$, \beq\PPP\po \mmi \wt{\Pi}(T\cdot)-W(\cdot)\mmi\geq
z+C_1\log T\pff\le C_2\exp\po-C_2z\pf.\label{kmtPoisson}\eeq Notice
that, as $k\rar \infty$~: 
\[\frac{(\nkkd a_{\nkkd})^{1/2}b_{\nkkd}\e_{\nkkd}}{\log(\nkkd
a_{\nkkd})} \rar \infty.\]  Thus, we have, ultimately as $k\rar
\infty$, \beq\PPP_k^{KMT} \le C_2\exp\poo-\frac{\e
C_3}{\sqrt{2}q_1p_2} (\nkkd a_{\nkkd})^{1/2}(\log_2
\nkkd)^{-1/6}\pff\label{etca}.\eeq Recalling the assumption
$na_n/(\log_2n) ^{7/3}\rar \infty$ we see that $\PPP_k^{KMT}$ is
sumable in $k$. Now, making use of assertion (\ref{rw}) of Fact \ref{propo: recouvrement Pbrownien} we have, for all large
$k$,
\begin{align}
\nonumber\PPP_k^{Tal}\;&=\;\PPP\poo W\notin
b_{\nkkd}\SSS_1+\frac{1}{q_1q_2}\e_{\nkkd}b_{\nkkd}\mcal{B}_0\pff \\
\nonumber \nonumber&\le \;\exp\pooo-\Big(\frac{\e}{q_1q_2}-\frac{L_1
(q_1q_2)^2}{2\e^2}\Big)(\log_2 \nkkd)^{1/3}-\log_2 \nkkd\pfff.
\end{align}
Now if $\e> (L_1/2)^{1/3}=: \e_1$ and if
$q_1,q_2>1$ are chose sufficiently small, then $\PPP_k^{Tal}$ is
sumable in $k$. By the Borel-Cantelli lemma, we see that for any
\(\e>\e_1\)  we have almost surely, for all large $n$,
\[g_n\in \SSS_1+\e_{\nkkd} \mcal{B}_0,\]
where $g_n:= (\nkkd
a_{\nkkd})^{-1/2}b_{\nkkd}^{-1}H_n(a_{\nkkd}\cdot ),\;\;n\in N_k$.
 To conclude the proof of Theorem \ref{theo: recproclocal}, it remains to control the distance between $a_n^{-1/2}b_n^{-1}\alp_n(a_n\cdot) $ and \(g_n\), which is the purpose of the following lemma.\lb
\begin{lem}\label{lem: kontroll}
We have almost surely~: 
\[\lsn \;(\log_2n)^{2/3}\Mmi \frac{\alp_n(a_n\cdot)}{(2a_n\log_2 n)^{1/2}}-g_n\Mmi =0.\]
\end{lem}
\textbf{Proof}~:
 Set
$\Gam_n:=1-(n/\nkkd)^{1/2}(a_n/a_{\nkkd})^{1/2}(\log_2n/\log_2
\nkkd)^{1/2}.$ The triangle inequality yields \begin{align} \Mmi
\frac{\alp_n(a_n\cdot)}{(2a_n\log_2 n)^{1/2}}-g_n\Mmi&\le\; \Mmi
\frac{\alp_n(a_n\cdot)}{(2a_n\log_2n)^{1/2}}\Gam_n\Mmi+\Mmi
\frac{H_n(a_n\cdot)-H_n(a_{\nkkd}\cdot)}{(2\nkkd a_{\nkkd}\log_2
a_{\nkkd})^{1/2}}\Mmi
\\ \nonumber &:=\;A_n+B_n.
\end{align}
Clearly we have, as $\kif$,
\[\max_{n\in N_k}\;(\log_2\nkkd)^{2/3}\;\Gam_n\le
\pooo 1-\sqrt{\frac{\nkd \log_2\nkd}{\nkkd \log_2 \nkkd}}\pfff
(\log_2\nkkd)^{2/3}\rar 0.\] Now, by applying (\ref{flilmason}) we
have almost surely \beq \lsn \Mmi
\frac{\alp_n(a_n\cdot)}{(2a_n\log_2
n)^{1/2}}\Mmi=1\label{m2norme}.\eeq Obviously (\ref{m2norme})
implies that, almost surely~: 
\[\limn (\log_2\nkkd)^{2/3}\max_{n\in  N_k}A_n=0.\]
 We now focus on controlling $B_n$. Set
$\rho_k:=a_{\nkd}/a_\nkkd$ and notice that
\begin{align}
\nonumber&\PPP\pooo\max_{n\in N_k}(\log_2 n)^{2/3}\Mmi
\frac{H_n(a_n\cdot)-H_n(a_{\nkkd}\cdot)}{(2\nkkd a_{\nkkd}\log_2
a_{\nkkd})^{1/2}}\Mmi \geq \e\pfff\\ \nonumber \le\; & \PPP\poooo
\max_{n\in N_k}\sup_{1\le\rho\le \rho_k,\;0\le t\le
1}(\log_2\nkkd)^{2/3}\Mid\frac{\alp_n(a_\nkkd \rho
t)-\alp_n(a_{\nkkd} t)}{(2a_{\nkkd} \log_2 \nkkd)^{1/2}}\Mid\geq
\e\pffff.  \end{align} Now consider the Banach space
$B\po[0,1]\times[0,2]\pf$ of all real bounded functions on
$[0,1]\times[0,2]$, endowed with the usual sup norm
$\norm_{[0,1]\times[0,2]}$. We shall now make use of the powerful
maximal inequality of Montgommery-Smith. For fixed $k\geq 1$, we
apply the just mentioned inequality to the finite sequence $(X_i)_{i\in N_k}$, with
$X_i(t,\rho):=\mathds{1}_{[t,\rho t]}(U_i)-\rho t,\;t\in[0,1],\;\rho\in
[1,\rho_k],\:\rho t \le 1$ and $X_i(t,\rho)=0$ elsewhere. Hence, by a combination of Theorem 1 and Corollary 3 in \cite{Montgom}, we have~: 
\begin{align} \nonumber&\PPP\pooo \max_{n\in N_k}\sup_{1\le\rho\le
\rho_k,\;0\le t\le 1}(\log_2\nkkd)^{2/3}\Mid\frac{\alp_n(a_\nkkd
\rho t)-\alp_n(a_{\nkkd} t)}{(2a_{\nkkd} \log_2
\nkkd)^{1/2}}\Mid\geq \e\pfff\\
\nono \le&9\PPP\pooo \sup_{1\le\rho\le \rho_k,\;0\le t\le
1}(\log_2\nkkd)^{2/3}\Mid\frac{\alp_{\nkkd}(a_\nkkd \rho
t)-\alp_{\nkkd}(a_{\nkkd}t)}{(2a_{\nkkd} \log_2
\nkkd)^{1/2}}\Mid\geq \e/30\pfff\\
 \nono\le\;& 18\;\PPP\pooo \mmi
\wt{\Pi}(\nkkd a_{\nkd}\cdot)-W(\nkkd a_{\nkd}\cdot)\mmi\geq
\frac{\e}{240}\frac{(2\nkkd a_\nkkd \log_2\nkkd)^{1/2}}{(\log_2 (
\nkkd))^{2/3}}\pfff\\  &+18\PPP\pooo\Mmi
\frac{W(\rho_k\cdot)-W(\cdot)}{(2 \log_2\nkkd)^{1/2}}\Mmi\geq
\frac{\e}{120(\log_2\nkkd)^{2/3}}\pfff\label{dalm}.
\end{align}
In the last expression (which is the combination of usual
poissonization techniques with the triangular inequality),
$\wt{\Pi}$ and $W$ denote respectively a centered Poisson process
and a Wiener process based on the same underlying probability
space. By the Koml\'os-Major-Tusnàdy construction (see \cite{KMT2}), $W$ can be
constructed to satisfy (\ref{kmtPoisson}). By making use of the
same arguments as those invoked to obtain (\ref{etca}), we
conclude that the first term in (\ref{dalm}) is sumable in
$k$. To control the second term in (\ref{dalm}), we shall
 make use of a well
known inequality (see, e.g., \cite{MSW}, p. 536), with
$a:=\rho_k-1$,
$\lab:=(\rho_k-1)^{-1/2}(\log_2\nkkd)^{-1/6}(\sqrt{2}\e/120)$
and $\dd:=1/2$, to get \begin{align}\nono &\PPP\poo\Mmi
\frac{W(\rho_k\cdot)-W(\cdot)}{\sqrt{2 \log_2\nkkd}}\Mmi\geq
\frac{\e}{120(\log_2\nkkd)^{2/3}}\pff\\
\nono\le &\frac{30720}{\sqrt{2}\e}(\rho_k-1)^{-1/2}((\log_2
\nkkd)^{1/6}\exp\poo-\frac{\e^2}{19200}(\rho_k-1)^{-1}(\log_2 \nkkd)^{-1/3}\pff. \end{align}
This expression is sumable in $k$, and hence $\max_{n\in
N_k}B_n\rar 0$ almost surely as $\kif$.$\Box$
\subsection{Proof of Theorem \ref{theo:
Chunglocal}} Recall that $\Cf$, $\Nf$, $\SSS_1^{bv}$ and
$\SSS_1^{liv}$ are defined in \S 2. The main tool to achieve the
proof of Theorem \ref{theo: Chunglocal} is the following inequality
(see Berthet \cite{Berthet2}), which sums up different results from Gorn and
Lifshits \cite{GornLif}, Berthet and Lifshits \cite{BerthetLif} and
Grill \cite{Grill1} (see also de Acosta \cite{Acosta1}).
\begin{ineg}\label{ineg: petiteboule} For any  $f\in \SSS_1^{BV}\cup\SSS_1^{LIV}$ and
$\dd>0$, there exist $\gam^+=\gam^+(\dd,f)>0$ and
$\gam-=\gam^-(\dd,f)>0$ such that for $T$ sufficiently large~: 
\[\PPP\pooo \Nf\poo\frac{T^2}{2}\pff\Mmi \frac{W}{T}-f\Mmi \le
(1+\dd)\Cf\pfff\geq
\exp\poo-\frac{T^2}{2}+\gam^+\frac{\Nf^2(T^2/2)}{T^2}\pff,\]
\[\PPP\pooo \Nf\poo\frac{T^2}{2}\pff\Mmi \frac{W}{T}-f\Mmi \le
(1-\dd)\Cf\pff\le
\exp\poo-\frac{T^2}{2}-\gam^-\frac{\Nf^2(T^2/2)}{T^2}\pfff.\]
\end{ineg}
Select $f\in \SSS_1^{BV}\cup\SSS_1^{LIV}$. We remind the two
following properties of $\Nf$ (see \cite{Tal4}), namely $
\ls_{L\rar \infty} L^{-1}\Nf(L)<\infty$ and $\li_{L\rar \infty}
L^{-2/3}\Nf(L)>0$. We shall first show that, almost surely~:
\[\lin\Nf(\log_2n ) \Mmi\frac{\alp_n(a_n
\cdot)}{(2a_n\log_2n)^{1/2}}-f\Mmi\geq \Cf.\]
 Let us fix $\e>0$. We start by applying poissonization techniques in combination with the Koml\'os-Major-Tusnàdy
 approximation.
\begin{align}\nono&\PPP\poo\Nf(\log_2n )\Mmi\frac{\alp_n(a_n
\cdot)}{(2a_n\log_2n)^{1/2}}-f\Mmi\le \Cf(1-2\e)\pff\\
\nono\le& 2\PPP\poo\Nf(\log_2n )\Mmi
\frac{W(na_n\cdot)}{(2na_n\log_2n)^{1/2}}-f\Mmi\le
\Cf(1-\e)\pff\\
\nono &+2\PPP\poo\Mmi W(na_n\cdot)-\wt{\Pi}(na_n\cdot)\Mmi \geq
\frac{\Cf\e(2na_n\log_2n)^{1/2}}{\Nf(\log_2n )}\pff. \end{align}
These two terms are sumable along the subsequence $\nkd$, the
second term being controlled by the Koml\'os-Major-Tusnàdy approximation while the
first one is controlled by Inequality  \ref{ineg: petiteboule}. Now
the control between $\nkd$ and $\nkkd$ follows the same line as in
Lemma \ref{lem: kontroll}. We omit details for sake of brevity.\lb
We now focus on showing that, almost surely,
\[\lin\Nf(\log_2n ) \Mmi\frac{\alp_n(a_n
\cdot)}{(2a_n\log_2n)^{1/2}}-f\Mmi\le \Cf.\] Set
$n_k:=k^{2k},\;\;v_k:=n_{k+1}-n_k$ and
\beq h_k:=\frac{\sqrt{n_{k+1}}\alp_{n_{k+1}}(a_{n_{k+1}}\cdot)-\sqrt{n_k}\alp_{n_k}(a_{n_{k+1}}\cdot)}{\sqrt{2v_ka_{n_{k+1}}\log_2(
v_k)}}.\label{hk}\eeq Notice that the $h_k$ are mutually independent processes,
and that each $h_k$ is distributed like
$(2a_{n_{k+1}}\log_2v_k)^{-1/2}\alp_{v_k}(a_{n_{k+1}}\cdot)$. We now
make use of the following "depoissonization" lemma. Recall that
$\wt{\Pi}(\cdot)$ denotes a centered standard Poisson process on
$[0,\infty)$.
\begin{lem}\label{lem: depoissonisation} Under assumptions (\ref{ww}) and (\ref{www}),
there exist two sumable positive  sequences $(c_k)_{k\geq
1},\;(b_k)_{k\geq 1}$ and an integer $k_0\geq 1$ such that, for
any set $A\subset B([0,1])$ that is measurable for both empirical
and Poisson processes and for all $k\geq k_0$,
\[\PPP\poo\wt{\Pi}(v_k a_{n_{k+1}}\cdot)\in A\pff-c_k-b_k\le
2\PPP\poo v_k^{1/2}\alp_{v_k}(a_{n_{k+1}}\cdot)\in A\pff.\]
\end{lem}
\textbf{Proof :}
Set $u_k:=(4\log_2( v_k)/n_{k+1}a_{n_{k+1}})^{1/2}$. By
assumption (\ref{www}) we have $u_k\rar 0$ as $\kif$. Now set
$\Pi(t):=\wt{\Pi}(t)+t,\;\in [0,1]$ and $R_{1,k}:=\Pi(v_k
a_{n_{k+1}})$, $R_{2,k}:=\Pi(v_k)-\Pi(v_k a_{n_{k+1}}).$ For fixed
$k$, $R_{1,k}$ and $R_{2,k}$ are independent random variables and
are distributed as Poisson variables with respective expectations
$v_k a_{n_{k+1}}$ and $v_k(1-a_{n_{k+1}})$. Let $A\subset B([0,1])$
be an arbitrary set that we assume to be measurable for $\wt{\Pi}$
and $\alp_n$. Define the following events~: 
\[E_k:=\aoo\wt{\Pi}(v_k a_{n_{k+1}}\cdot)\in A\aff,\;k\geq 1.\]
We have, for any integer $k\geq 1$,
\begin{align}
\nonumber\PPP(E_k) \le\;& \PPP\poo E_k\cap R_{1,k}\in
[(1-u_k)v_ka_{n_{k+1}},(1+u_k)v_ka_{n_{k+1}}]\pff \\ \nonumber
\nonumber &+\PPP\poo R_{1,k}<(1-u_k)v_ka_{n_{k+1}}\pff+\PPP\poo
R_{1,k}>(1+u_k)v_ka_{n_{k+1}}\pff.
\end{align}
Denote by $c_k$ and  $b_k$ the two last terms of the RHS of the
preceding inequality. We shall show that these two sequences have
finite sums. Making use of Chernoff's inequality, we have~: 
\[c_k\le \exp\poo-v_ka_{n_{k+1}}\poo(1+u_k)\log(1+u_k)-u_k\pff\pff.\]
Since $(1+u)\log(1+u)-u\sim \frac{u^2}{2}$ as $u\rar 0$, it
follows that for all large $k$, \begin{align} \nono c_k\le&
\exp\poo
-v_ka_{n_{k+1}}\frac{u_k^2}{2}\pff\\
\nono=&\exp\poo-2\log_2n_k\pff. \end{align} We make use of a
similar method to show that $(b_k)_{k\geq 1}$ is sumable. It
remains to show that, for all $k\geq k_0$ (with $k_0$ independent
of $A$), we have
\[\PPP\poo E_k\cap R_{1,k}\in
[(1-u_k)v_ka_{n_{k+1}},(1+u_k)v_ka_{n_{k+1}}]\pff\le 2\PPP\poo
E_k\mid \Pi(v_k)=v_k\pff.\] Now set \beq K_k:=\inf \aooo
\frac{\PPP(R_{2,k}=v_k-j)}{\PPP(\Pi(v_k)=v_k)}j\in[(1-u_k)v_ka_{n_{k+1}},(1+u_k)v_ka_{n_{k+1}}]\afff.\label{wfm}\eeq
We have
\begin{align}
\nonumber&\PPP\poo E_k\cap R_{1,k}\in
[(1-u_k)v_ka_{n_{k+1}},(1+u_k)v_ka_{n_{k+1}}]\pff \\ \nonumber
\nonumber\le\;&\sli_{j=[(1-u_k)v_ka_{n_{k+1}}]}^{[(1+u_k)v_ka_{n_{k+1}}]+1}\PPP(E_k\cap
R_{1,k}=j) \\ \nonumber \nonumber\le \;& K_k^{-1}
\sli_{j=[(1-u_k)v_ka_{n_{k+1}}]}^{[(1+u_k)v_ka_{n_{k+1}}]+1}\PPP(E_k\cap
R_{1,k}=j)\frac{\PPP(R_{2,k}=v_k-j)}{\PPP(\Pi(v_k)=v_k)} \\
\nonumber
\nonumber=\;&K_k^{-1}\sli_{j=[(1-u_k)v_ka_{n_{k+1}}]}^{[(1+u_k)v_ka_{n_{k+1}}]+1}\frac{\PPP(E_k\cap
R_{1,k}=j\cap R_{2,k}=v_k-j)}{\PPP(\Pi(v_k)=v_k)} \\ \nonumber
\nonumber\le\;& K_k^{-1}\PPP(E_k\mid \Pi(v_k)=v_k).
\end{align}
Hence, it suffices to show that $K_k\rar 1$. For clarity of notations, set  $v'_k:=v_k(1-a_{n_{k+1}})$, recalling
that $R_{2,k}$ is a Poisson variable with parameter $v'_k$. Setting $l=v_k-j$ in
(\ref{wfm}) we have, as $\kif$,
\[K_k=\inf \aooo \frac{\PPP(R_{2,k}=l)}{\PPP(\Pi(v_k)=v_k)}, l\in
[v'_k-v_ku_k a_{n_{k+1}},v'_k+v_ku_ka_{n_{k+1}}]\afff.\]
Now, by
Stirling's formula, we have $\PPP(\Pi(v_k)=v_k)\sim (2\pi
v_k)^{-1/2}$ as $\kif$. A routine study of the finite sequence
\[\poo\PPP(R_{2,k})=l,\;l\in
[v'_k-v_ku_k a_{n_{k+1}},v'_k+v_ku_ka_{n_{k+1}}]\pff\]
shows that
\begin{align}
\nono&\PPP(\Pi(v_k)=v_k) K_k=\min (\PPP_{1,k},\PPP_{2,k})\text{ , where}\\
\nono&\PPP_{1,k}:=\PPP(R_{2,k}=[v'_k-v_ka_{n_{k+1}}u_k])\text{, and}\\
\nono&\PPP_{2,k}:=\PPP(R_{2,k}=[v'_k+v_ka_{n_{k+1}}u_k]+1).
\end{align}
We set $u'_k=a_{n_{k+1}}u_k v_k/v'_k\sim u_k a_{n_{k+1}}$.
Stirling's formula yields, ultimately as $\kif$,
\begin{align}
\nonumber\PPP_{1,k}&=\;\frac{{{v_k}'}^{[{v_k}'-v_ka_{n_{k+1}}u_k]}}{[{v_k}'-v_ka_{n_{k+1}}u_k]!}\exp(-{v_k}')
\\ \nonumber
\nonumber&\sim\; {\poo
\frac{{v_k}'}{{v_k}'(1-{u_k}')}\pff}^{{v_k}'(1-{u_k}')}\frac
{\exp(-{v_k}')}{\exp(-{v_k}'+{v_k'}{u_k'})}\sqrt{2\pi v'_k} \\
\nonumber \nonumber&\sim\;\sqrt{2\pi
v_k}(1-{u_k}')^{-{v_k'}(1-{u_k}')}\exp(-{v_k}'{u_k}') \\ \nonumber
\nonumber&=\;\sqrt{2\pi
v_k}\exp\poo-{v_k}'\po(1-{u_k}')\log(1-{u_k})+{u_k}'\pf\pff.
\end{align}
Moreover, since $(1-\e)\log(1-\e)+\e\sim\frac{\e^2}{2}$ as $\e\rar
0$, we have, for all large $k$, $$\exp\poo-2a_{n_{k+1}}\log_2(
n_{k+1})\pff\le
\exp\poo-{v_k'}\po(1-{u_k}')\log(1-{u_k})+{u_k}'\pf\pff.$$ By
assumption (\ref{ww}) we have $a_n\log_2n \rar 0$, which ensures
that $\PPP_{1,k}\sim\sqrt{2\pi v_k}$. The control of $\PPP_{2,k}$
is very similar. This achieves the proof of Lemma \ref{lem:
depoissonisation}.$\Box$\lb
We now apply the preceding lemma in
conjunction with the Koml\'os-Major-Tusnàdy approximation. Let $W$ be a Wiener
process constructed on the same underlying probability space as
$\wt{\Pi}$. For an arbitrary $\dd>0$, we have (recall tha $h_k$ has been defined in (\ref{hk}))
\begin{align}
\nonumber&\PPP\poo\Nf (\log_2( v_k))\mmi h_k-f\mmi \le
(1+2\dd)\Cf\pff \\ \nonumber \nonumber\geq\;& \frac{1}{2}\PPP\pooo\Nf
(\log_2( v_k))\Mmi
\frac{\wt{\Pi}(v_ka_{n_{k+1}}\cdot)}{(2v_ka_{n_{k+1}}\log_2v_k)^{1/2}}-f\Mmi
\le (1+2\dd)\Cf\pfff-\frac{1}{2}c_k-\frac{1}{2}b_k \\ \nonumber
\nonumber\geq\;& -\frac{1}{2}\PPP\pooo\Mmi \wt{\Pi} (
v_ka_{n_{k+1}}\cdot)-W(v_ka_{n_{k+1}}\cdot)\Mmi\geq \frac{\dd
\Cf(2v_ka_{n_{k+1}}\log_2v_k)^{1/2}}{\Nf(\log_2v_k)}\pfff \\
\nonumber \nonumber&+\frac{1}{2}\PPP\pooo \Nf(\log_2v_k)\mmi
\frac{W}{(2\log_2v_k)^{1/2}}-f\mmi\le \Cf(1+\dd)\pfff-\frac{1}{2}c_k-\frac{1}{2}b_k \\
\nonumber
\nonumber=:\;&-d_k+\frac{1}{2}\PPP_k-\frac{1}{2}c_k-\frac{1}{2}b_k.
\end{align}
Since $v_k\sim n_{k+1}$ it is easy to conclude that $d_k$ is
sumable in $k$, by making use of the strong approximation (see
\cite{KMT2}). Hence, making use of Inequality \ref{ineg:
petiteboule}, we have asymptotically
\[\PPP_k\geq \exp\poo-\log_2v_k\pff.\]
But $\log_2v_k=\log (k+1) +\log_2 k+o\po k^{-2}(\log k)^{-1}\pf$ and
hence
\[\sum_{k\geq 1} \PPP\poo\Nf (\log_2( v_k))\mmi h_k-f\mmi \le
(1+2\dd)\Cf\pff=\infty.\] Applying the second half of the Borel
Cantelli lemma, we deduce that, almost surely~: 
\[\liminf_{k\rar\infty}\Nf(\log_2v_k)\mmi h_k-f\mmi \le \Cf.\]
To conclude the proof, it is enough to show  that, almost surely
(recall that $\lim_{L\rar \infty} L^{-1}\Nf(L)>0$),
\[\lim_{k\rar \infty} (\log_2 n_{k+1})\;\Mmi
h_k-\frac{\alp_{n_{k+1}}(a_{n_{k+1}}\cdot)}{(2n_{k+1}a_{n_{k+1}}\log_2(
n_{k+1}))^{1/2}}\Mmi=0.\] Routine algebra shows that
\begin{align}
\nonumber&(\log_2 n_{k+1})\Mmi
h_k-\frac{\alp_{n_{k+1}}(a_{n_{k+1}}\cdot)}{(2a_{n_{k+1}}\log_2(
n_{k+1}))^{1/2}}\Mmi \\ \nonumber \nonumber\le\;&(\log_2
n_{k+1})\poo\po\frac{n_{k+1}\log_2
n_{k+1}}{v_k\log_2(v_k)}\pf^{1/2}-1\pff\Mmi\frac{\alp_{n_{k+1}}(a_{n_{k+1}}\cdot)}{(2a_{n_{k+1}}\log_2
n_{k+1})^{1/2}}\Mmi\\
\nonumber&+(\log_2 n_{k+1})\Mmi
\frac{n_k^{1/2}\;\alp_{n_k}(a_{n_{k+1}}\cdot)}{(2v_k
a_{n_{k+1}}\log_2n_{k+1})^{1/2}}\Mmi
\\ \nonumber :=\;& A_k+B_k.
\end{align}
Applying theorem of Mason (\ref{flilmason}) we get $A_k \rar 0$
almost surely as $k\rar \infty$. We now apply Doob's inequality
for positive submartingales to obtain
\begin{align}
\nonumber&\PPP\poo (\log_2 n_{k+1})\Mmi
n_k\;\frac{\alp_{n_k}(a_{n_{k+1}}\cdot)}{((2v_k
a_{n_{k+1}}\log_2n_{k+1}))^{1/2}}\Mmi\geq \e\pff \\ \nonumber
\nonumber=\;&\PPP\poo \sup_{0\le t\le a_{n_{k+1}}}\Mid
\frac{\alp_n(t)}{1-t}\Mid\geq \frac{\e(2v_k
a_{n_{k+1}}\log_2n_{k+1}))^{1/2}}{n_k^{1/2}\log_2n_{k+1}}\pff \\
\nonumber \nonumber \le\;&\frac{1}{2\e^2}(1-a_{n_{k+1}}) \log_2
n_{k+1} \frac{n_k}{v_k}.
\end{align} Since $ n_k/v_k\sim
1/e^2k^2$ as $\kif$, we conclude the proof of the lower bound in
Theorem \ref{theo: Chunglocal} with the Borel-Cantelli lemma.
$\Box$\lb
\textbf{Acknowledgements}: The author would like to thank P.
Berthet and E. Shmileva for many fruitful conversations and
advices on the subject.

\end{document}